\documentclass[12pt]{iopart}
\expandafter\let\csname equation*\endcsname\relax
\expandafter\let\csname endequation*\endcsname\relax
\usepackage{mathtools}
\usepackage{graphicx}
\usepackage{float}
\usepackage{amsmath}
\usepackage{geometry}
\usepackage{bm}
\usepackage{caption}
\usepackage{subcaption}
\usepackage{tabularx}
\usepackage{array}
\usepackage{mathrsfs}
\usepackage{bbm}
\usepackage{epstopdf}

\usepackage{algorithm}
\usepackage{algpseudocode}

\algdef{SE}[DOWHILE]{Do}{doWhile}{\algorithmicdo}[1]{\algorithmicwhile\ #1}

\renewcommand\theequation{\arabic{section}.\arabic{equation}}
\numberwithin{equation}{section}
\newtheorem{thm}{Theorem}[section]
\newtheorem{lemma}{Lemma}[section]
\newtheorem{example}{Example}[section]

\begin{document}

\title[Heterogeneous and anisotropic elastic parameter estimation using DMOL]{Heterogeneous and anisotropic elastic parameter estimation using a novel semi-analytical forward solver}

\author[X. Zhu, and Z. Huang]{Xiaopeng Zhu\textsuperscript{1,$*$} and  Zhongyi Huang\textsuperscript{2}}

\address{\textsuperscript{1} College of Science, Liaoning Technical University, Liaoning, China}
\address{\textsuperscript{2} Department of Mathematical Sciences, Tsinghua University, Beijing, China}
\ead{zxp17@tsinghua.org.cn, zhongyih@mail.tsinghua.edu.cn}
\vspace{10pt}
\begin{indented}
\item[]April 2025 
\end{indented}

\begin{abstract}
An efficient procedure using a novel semi-analytical forward solver for identifying heterogeneous and anisotropic elastic parameters from only one full-field measurement is proposed and explored. We formulate the inverse problem as an special energy functional minimization with total variation(TV) regularization. The minimization problem is solved by Adam algorithm, which only requires solving one forward problem and no adjoint problem in each iteration. In order to deal with the irregularity of the elastic regions, the anisotropy and heterogeneity of parameters and potential singularities in forward-modeled issues, a novel semi-analytical forward solver named the direct method of lines is proposed, which discretizes angular variable while preserving analytical solutions along remaining coordinates. To validate the efficacy of our procedure, a series of numerical experiments are implemented subsequently, achieving reliable performance in both forward modeling and the six elastic arguments reconstruction scenarios. 

\vspace{2pc}

\end{abstract}

%
%
%
%
%

\section{Introduction} \label{sec:intro}
The estimation of material parameters involves identifying unknown coefficient (or parameter) functions based on imperfectly measured system responses. Solving such problems numerically holds substantial practical importance across multiple engineering disciplines, including seismic and biomechanical imaging, material property analysis, and structural integrity assessment. Nevertheless, conventional numerical methods frequently struggle due to the intrinsic ill-posedness and non-unique solutions characteristic of inverse problems. Consequently, numerous estimation approaches have been developed and documented in the literature. For elastic parameter estimation, detailed overviews of various identification methods are available in \cite{Avril2008,GBao2019,Bonnet2005,Fan2011}. Although techniques for identifying isotropic elastic parameters are relatively mature \cite{ Albocher2014,Banerjee2009,Banerjee2013,Jadamba2014,Oberai2003,Worden2007}, estimating anisotropic constitutive parameters from displacement fields, mode shapes, or strain measurements remains challenging. This difficulty arises due to the higher complexity and greater number of material variables involved, and a key challenge in anisotropic parameter estimation lies in acquiring measurements with sufficient sensitivity to enable accurate reconstruction. Determining the minimum number of physical tests required to fully characterize anisotropic properties remains an open research question. Furthermore, when material heterogeneity is taken into account, the need for multiple tests combined with inherent strain field variations introduces additional complexity, fundamentally altering the nature of the identification problem.

Significant progress in full-field measurement methodologies has facilitated diverse strategies for characterizing homogeneous anisotropic material properties. Initial studies by Bruno et al.\cite{Bruno2002,Bruno2005} pioneered experimental techniques employing full-field measurements under bending stress conditions to analyze both isotropic and anisotropic plate behavior. Subsequent methodological developments include Genovese's work\cite{Genovese2004}, which integrated speckle interferometry measurements with analytical models to extract in-plane elastic constants of fiber-reinforced polymer composites. In addressing orthotropic material systems, researchers have developed more efficient characterization approaches: Lecompte and colleagues \cite{Lecompte2007} established a framework combining finite element analysis with single-specimen biaxial testing, while alternative methodologies \cite{Molimard2005} based on open-hole tensile configurations were concurrently developed for practical stiffness evaluation.

Identifying heterogeneous and anisotropic elastic parameters remains significantly reduced explored compared to homogeneous cases, primarily due to the substantial increase in inverse unknowns within an inherently ill-posed optimization framework. Current approaches typically employ functionals under least-squares (L2) minimization leveraging gradient-driven optimization. Notable contributions include Liu et al.\cite{Liu2005}, who determined full anisotropic constants in heterogeneous media with known interior and exterior subdomain boundaries using computed tomography data, requiring complete boundary displacement data and partial external force information.  More recently, Genovese et al.\cite{Genovese2014} developed an inverse characterization approach for heterogeneous gallbladder tissue elasticity using digital image correlation, while \cite{Raghupathy2010} introduced a direct identification technique. Using an L2 minimization framework which augmented by semi-norm regularization, the work in \cite{Shore2011} effectively determined transverse anisotropic arguments of cancellous bone. Nevertheless, L2-based methods face notable limitations, including strong dependence on initial guesses and slow convergence with quasi-Newton methods.  Then in \cite{Guchhait2016}, the anisotropic linear elastic parameters were reconstructed by Guchhait et al. using a constitutive equation error minimization framework.

More recent approach for elastic parameter estimation has emerged, formulating the inverse problem as an energy functional minimization with total variation (TV) regularization \cite{part1,part2,ident-lame}. The minimization problem is usually solved by gradient-based optimization algorithm such as the Adam algorithm \cite{Kingma2015}. A heterogeneous and linear anisotropic elastic problem in an irregular domain is proposed in each iteration.  A significant challenge in solving the forward problem arises from stress singularities that typically develop at three critical locations: (1) material interface junctions, (2) geometric edges, and (3) crack tips. These singularities fundamentally complicate stress field analysis and numerical solution convergence.  From \cite{gonzalez2013efficient} we know that the singular part of displacement $u$ is $\mathcal{O}(r^{b_j})$ with $\textbf{Re}(b_j)\in(0,1)$ in the polar coordinate. Thus  conventional numerical approaches such as finite element methods (FEM) exhibit significantly diminished convergence rates\cite{chang2007singular,gonzalez2013efficient}. To improve the numerical result, some more efficient methods need to be developed.

Among existing approaches, the direct method of lines has emerged as particularly effective semi-discrete technique for addressing singular elliptic problems. As an advanced development of the classical method of lines \cite{schiesser2012numerical,xanthis1991method}, this approach offers two key advantages: (1) it requires no a priori knowledge of singularity characteristics, and (2) it naturally captures solution singularities. Originally developed for handling elliptic problems with singularities in polygonal domains \cite{H-H-1999,com-material} and boundary value problems on unbounded domains \cite{H-W-1996,H-W-1999}, the method has demonstrated remarkable versatility. Most recently in star-shaped domains, its applicability was extended in \cite{W-H-2018} to solve Laplace's equation with singularities and to solve compressible\cite{zhu2023}  and nearly-incompressible \cite{zhu2024} elasticity problems of composite materials.

In this paper, an efficient computational framework for identifying heterogeneous and  anisotropic elastic parameters using single full-field measurements will be presented.  Its forward solver is a novel semi-analytical scheme named the direct method of lines, offering distinct advantages for handling material anisotropy, parameter heterogeneity, irregular domain geometries, and potential solution singularities. Our inverse formulation minimizes a specialized energy functional under a least-squares criterion with total variation regularization, implemented via the Adam optimization algorithm. This configuration provides computational efficiency by requiring only one forward solution per iteration, eliminating the need for adjoint problem solutions. Numerical experiments will demonstrate the method's effectiveness in both forward simulation and parameter identification tasks, validating its accuracy across various test cases. 

The remainder of this paper is structured as follows. Section \ref{sec:inverse} presents our novel computational framework for identifying heterogeneous and anisotropic elastic parameters via only one full-field displacement measurement. Section \ref{sec:GDMOL} details the implementation of the semi-analytical scheme as an efficient forward solver, and an optimal error estimate will be achieved. Comprehensive numerical experiments in Section \ref{sec:example} validate the method's accuracy and robustness through both forward and inverse problem solutions. Finally, Section \ref{sec:conclusion} summarizes key findings and discusses potential extensions of this work.

\section{Heterogeneous and anisotropic elastic parameter estimation problem }\label{sec:inverse}

\begin{figure}[H]
\centering
\includegraphics[scale=0.3]{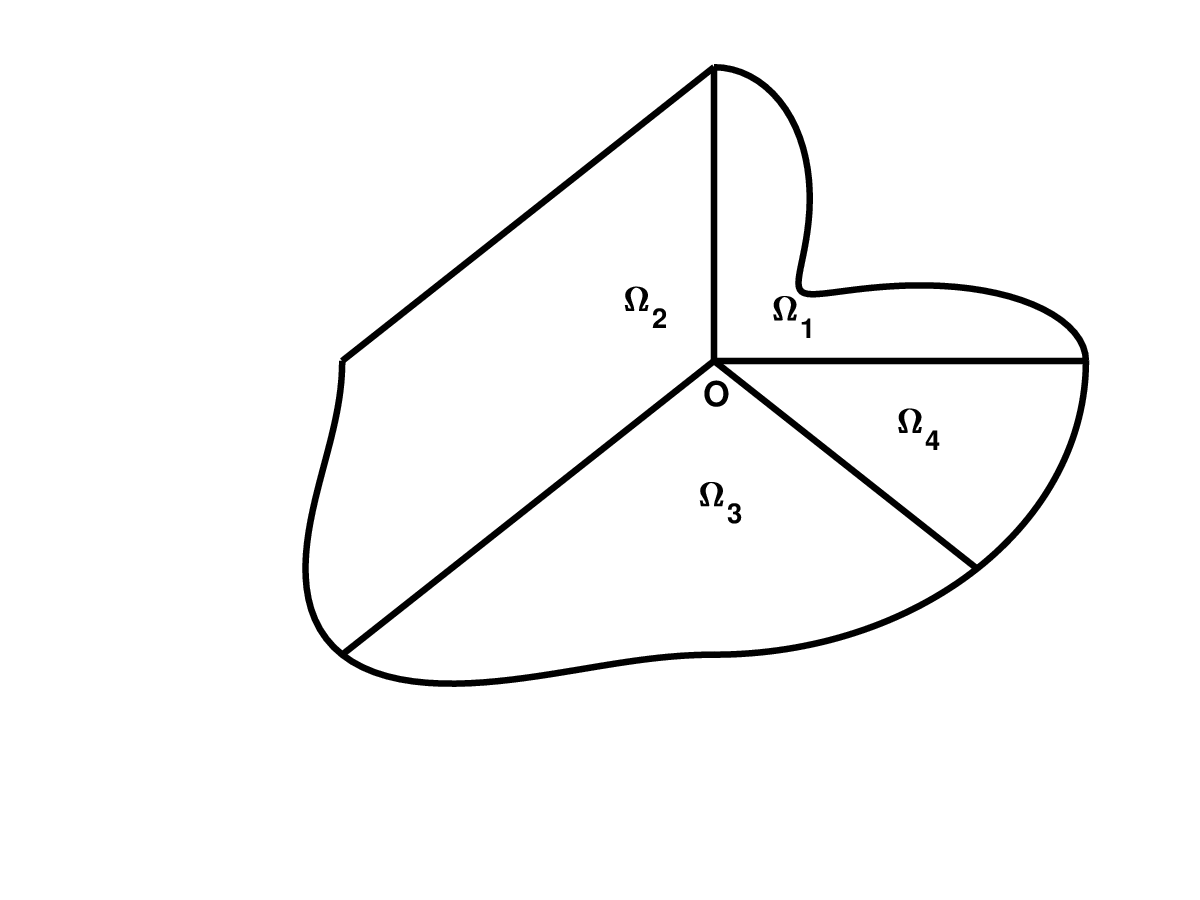}\vspace{-5mm}
\caption{Heterogeneous material systems within star-configured regions}
\label{fig:domain0}
\end{figure}

We are interested in the inverse linear elasticity problem to anisotropic compressible composite materials in star-shaped regions, which is also a heterogeneous and anisotropic elastic parameter estimation problem.
We begin by considering a star-shaped domain $\overline{\Omega} = \bigcup_{k = 1}^K \overline{\Omega}_k \subset R^2$, where $\Omega$ represents the entire composite material consisting of $K$ distinct material phases, and $\Omega_k$ denotes the $k$-th material phase. The boundary $\Gamma = \partial \Omega$ is assumed to be star-configured relative to $O$ and $\tilde{r}(\phi)$ provides a piecewise C¹ parametric description of the boundary, satisfying $\tilde{r}(0) = \tilde{r}(2\pi)$ with $\forall \phi\in[0,2\pi],0<r_0\leq\tilde{r}(\phi)$. The geometric property ensures the proper definition of the curvilinear coordinate system established with \eqref{coordinate}. No generality is lost when presuming different material surfaces intersect at $O$, with each surface $\overline{\Omega}_{k-1} \cap \overline{\Omega}_k$ given by a radial straight curve $L_k = \{(r, \theta) |0 \leq r \leq \tilde{r}(\theta_k), \theta=\theta_k \}$ with $k = 1, \ldots, K$ (see Figure \ref{fig:domain0}). Within this domain, we study the Navier equations under Dirichlet boundary conditions, which govern the equilibrium state of an anisotropic elastic medium:
\begin{equation} \label{composite-eq}
\begin{aligned}
-\nabla\cdot \sigma^k & = 0, \quad \text{in } \Omega_k,\; k = 1, \ldots, K\\
u^k & = f^k, \quad \text{on } \Gamma_k=\partial \Omega_k \bigcap \partial \Omega,\\
u^{k-1} = u^{k}, \sigma^{k-1}&\cdot {n}_k = \sigma^{k}\cdot {n}_k, \quad \text{on } L_k,
\end{aligned}
\end{equation}
with
\begin{equation} \label{definition of stress tensor}
\begin{aligned}
\varepsilon(u^k) &= \frac{1}{2} \left( \nabla u^k + (\nabla u^k)^T \right),
\left(
\begin{matrix}
\sigma^k_{11} \\
\sigma^k_{22} \\
\sigma^k_{12}
\end{matrix}\right)& = \left(
\begin{matrix}
a^k_{11} & a^k_{12} & a^k_{13}\\
a^k_{21} & a^k_{22} & a^k_{23}\\
a^k_{31} & a^k_{32} & a^k_{33}
\end{matrix}\right)\left(
\begin{matrix}
\varepsilon^k_{11} \\
\varepsilon^k_{22} \\
2\varepsilon^k_{12}
\end{matrix}\right),
\end{aligned}
\end{equation}
$f^k = (f^k_1, f^k_2)^T=f|_{\Gamma_k}$ with $f$ denotes the prescribed boundary condition, $u^k = (u_1^k, u_2^k)^T=u|_{\Omega_k}$ with $u$ denotes deformation field, $\sigma^k=\sigma|_{\Omega_k}$ with $\sigma$ represents tension matrix and the unit normal vector ${n}_k = (- \sin(\theta_k), \cos(\theta_k))^T$ to the interface $L_k$. The anisotropic coefficient tensor $a_{ij}^k$ satisfies the symmetric positive definiteness condition and is piecewise constant within each subdomain $\Omega_k$. For the orthotropic case, $a_{13}^k=a_{23}^k=0$. In the isotropic scenario, the coefficients reduce to $a_{11}^k=a_{22}^k=2\mu^k+\lambda^k,a_{33}=\mu^k,a_{12}=\lambda^k,a_{13}=a_{23}=0$, where $\lambda^k$ and $\mu^k$ are Lam\'e coefficients for material $k$.

For the system \eqref{composite-eq}, the forward problem consists of computing the displacement field $u$ given the anisotropic material coefficients $a$. Conversely, the inverse problem for \eqref{composite-eq} involves reconstructing the anisotropic coefficients $a$ from measured displacement data $u$.

To formulate the inverse elasticity problem for anisotropic composite materials as a regularized minimization problem. Let $u[a]$ denote the solution to \eqref{composite-eq} corresponding to the anisotropic coefficient tensor $a$. Given our focus on scenarios involving discontinuous anisotropic coefficients across material interfaces, we make the following key assumption
\begin{equation}
a\in\Lambda := \left\{ a\text{ is a self-adjoint matrix with strictly positive spectrum } | \text{TV}(a) < \infty \right\},
\end{equation}
where $\text{TV}(\cdot)$ denotes the TV semi-norm \cite{vogel2002computational}. In the heterogeneous and anisotropic elastic parameter estimation problem with full-field displacement measurements. Considering noise-perturbed measurement data $z(x, y)$ of the true displacement field $u[a^{\star}](x, y)$, where $a^{\star} \in \Lambda$ represents the unknown true anisotropic parameter tensor and $(x, y) \in \Omega$ spans the spatial domain.

To reconstruct $a^{\star}$, set
$$
\begin{aligned}
&J_0(a) = \frac{1}{2} \iint_{\Omega} \bigg( a_{11}(\frac{\partial(u_1[a]-z_1)}{\partial x})^2+a_{22}(\frac{\partial(u_2[a]-z_2)}{\partial y})^2+a_{33}(\frac{\partial(u_1[a]-z_1)}{\partial y}+\\
&\frac{\partial(u_2[a]-z_2)}{\partial x})^2 +2a_{12}\frac{\partial(u_1[a]-z_1)}{\partial x}\frac{\partial(u_2[a]-z_2)}{\partial y}+ 2a_{13}\frac{\partial(u_1[a]-z_1)}{\partial x}(\frac{\partial(u_2[a]-z_2)}{\partial x}\\
&+\frac{\partial(u_1[a]-z_1)}{\partial y})
+2a_{23}\frac{\partial(u_2[a]-z_2)}{\partial y}(\frac{\partial(u_2[a]-z_2)}{\partial x}+\frac{\partial(u_1[a]-z_1)}{\partial y}) \bigg) dxdy.
\end{aligned}
$$
Since computing the functional gradient of $J_0(a)$ does not require evaluating partial derivatives of the displacement field $u$ with respect to the anisotropic parameters $a$, this approach offers significant computational advantages. Each gradient evaluation requires solving just one forward elastic problem (Navier's equations), avoiding the need for additional adjoint-state solutions.
Consider
\begin{equation} \label{minpro}
\begin{aligned}
\min_{a \in \Lambda}J(a):= J_0(a) + \eta  \text{TV}(a) ,
\end{aligned}
\end{equation}
where $\eta > 0$ is the regularization parameter.
\begin{thm}\label{theory 1}
The minimization problem \eqref{minpro} exists at least one solution $a^0\in \Lambda$. 
\end{thm}
This result follows from Theorem 4.3 in \cite{part1} with our anisotropic elasticity setting.

When considering the problem in the curvilinear coordinate system \eqref{coordinate}, we assume that the true coefficient $a^{\star}$  is constant along radial ($\rho$) directions for each fixed angular $\phi$. In the angular ($\phi$) direction, $a^{\star}$ is piecewise constant with material interfaces at predetermined angles. Let $h = \max \limits_{0 \leq i< m} |\varphi_{i+1} - \varphi_i|$ and $0 = \varphi_0 < \varphi_1 < \ldots < \varphi_m = 2\pi$,
\begin{equation} \label{discrete admissible set}
\Lambda_h = \left\{ a_h(\rho, \phi) \in \Lambda \Bigg{|}
\begin{aligned}
& a_h \text{ maintains constant values almost surely} \\
& \text{in } (-\infty,0] \times [\varphi_i, \varphi_{i+1}) 
\end{aligned}
 \right\}.
\end{equation}
 Approximating the continuous problem \eqref{minpro} through the following discretization scheme
\begin{equation} \label{num-minpro}
\begin{aligned}
\min_{a_h \in \Lambda_h} J(a_h)  = J_0(a_h)+ \eta  \text{TV}(a_h).
\end{aligned}
\end{equation}
\begin{thm} \label{theory 2}
The discrete minimization problem \eqref{num-minpro} admits a solution $a_h^0 \in \Lambda_h$.
\end{thm}
which follows from Theorem \ref{theory 1} by observing the discrete space $\Lambda_h$ remains $L^1$-closed within the continuous space $\Lambda$.

Then we use the Adam algorithm \cite{Kingma2015} to solve the minimization problem \eqref{num-minpro}. Assume
\begin{equation} \label{discrete Lame coefficients}
a^h_{ij}(\rho, \phi) = \mathbbm{1}_{(-\infty, 0]}(\rho) \sum_{t = 1}^m a_{ij}^t \mathbbm{1}_{[\varphi_{t - 1},\varphi_t)}(\phi), i,j=1,2,3.
\end{equation}
And then
$$
\text{TV}(a_{ij}^h) = \sum_{t=1}^{m} \tilde{r}(\varphi_t) |a_{ij}^{t+1}-a_{ij}^t|,
$$
where $a_{ij}^{m + 1} = a_{ij}^1$. For the implementation of the Adam algorithm to solve \eqref{num-minpro}, the optimization procedure necessitates evaluating $\nabla J$, which in turn requires calculating $\nabla[\text{TV}(a_{ij}^h)]$. For $t=1, \ldots, m$,
\[
\frac {\partial TV(a_{ij}^h)}{\partial a_{ij}^t}=sgn(a_{ij}^t-a_{ij}^{t-1})\tilde{r}(\theta_{j-1})+sgn(a_{ij}^{t+1}-a_{ij}^t)\tilde{r}(\theta_{j}),
\]
where $a_{ij}^0 = a_{ij}^m, a_{ij}^{m + 1} = a_{ij}^1$ and $sgn(x)$ is the sign function.
Using the chain rule, we obtain
$$
\begin{aligned}
\frac {\partial J_0}{\partial a_{11}^t} &  = \int_{-\infty}^0 \int_{\varphi_{t-1}}^{\varphi_t} \frac12 ((\frac{\partial z_1}{\partial x})^2-(\frac{\partial u_1[a_h]}{\partial x})^2) e^{2\rho} \tilde{r}(\phi)^2 d\phi d\rho, \\
\frac {\partial J_0}{\partial a_{22}^t} &  = \int_{-\infty}^0 \int_{\varphi_{t-1}}^{\varphi_t} \frac12 ((\frac{\partial z_2}{\partial y})^2-(\frac{\partial u_2[a_h]}{\partial y})^2) e^{2\rho} \tilde{r}(\phi)^2 d\phi d\rho, \\
\frac {\partial J_0}{\partial a_{33}^t} &  = \int_{-\infty}^0 \int_{\varphi_{t-1}}^{\varphi_t} \frac12 ((\frac{\partial z_1}{\partial y}+\frac{\partial z_2}{\partial x})^2-(\frac{\partial u_1[a_h]}{\partial y}+\frac{\partial u_2[a_h]}{\partial x})^2) e^{2\rho} \tilde{r}(\phi)^2 d\phi d\rho, \\
\frac {\partial J_0}{\partial a_{12}^t} &  = \int_{-\infty}^0 \int_{\varphi_{t-1}}^{\varphi_t} (\frac{\partial z_1}{\partial x}\frac{\partial z_2}{\partial y}-\frac{\partial u_1[a_h]}{\partial x}\frac{\partial u_2[a_h]}{\partial y}) e^{2\rho} \tilde{r}(\phi)^2 d\phi d\rho, \\
\frac {\partial J_0}{\partial a_{13}^t} &  = \int_{-\infty}^0 \int_{\varphi_{t-1}}^{\varphi_t} (\frac{\partial z_1}{\partial x}(\frac{\partial z_2}{\partial x}+\frac{\partial z_1}{\partial y})-\frac{\partial u_1[a_h]}{\partial x}(\frac{\partial u_2[a_h]}{\partial x}+\frac{\partial u_1[a_h]}{\partial y})) e^{2\rho} \tilde{r}(\phi)^2 d\phi d\rho, \\
\frac {\partial J_0}{\partial a_{23}^t} &  = \int_{-\infty}^0 \int_{\varphi_{t-1}}^{\varphi_t} (\frac{\partial z_2}{\partial y}(\frac{\partial z_2}{\partial x}+\frac{\partial z_1}{\partial y})-\frac{\partial u_2[a_h]}{\partial y}(\frac{\partial u_2[a_h]}{\partial x}+\frac{\partial u_1[a_h]}{\partial y})) e^{2\rho} \tilde{r}(\phi)^2 d\phi d\rho.
\end{aligned}
$$
Hence for $t=1, \ldots, m$
\begin{equation} \label{gradient of J}
\frac {\partial J}{\partial a_{ij}^t}=\frac {\partial J_0}{\partial a_{ij}^t} + \eta\frac {\partial \text{TV}(a_{ij}^h)}{\partial a_{ij}^t}.
\end{equation}

From the preceding analysis, evaluating the objective functional $J(a_h)$ and its gradient requires computing the displacement field $u[a_h]$ for discrete parameters $a_h \in \Lambda_h$. Two critical challenges arise in this computation: (1) $u[a_h]$ potentially contains singular stress behavior at material interface junctions and the geometric singularity point (where all interfaces $L_k$ meet), (2)the star-configured regions $\Omega$ displays inherent geometrical irregularities. Hence the forward solver is a novel semi-analytical scheme named the direct method of lines, chosen for its ability to handle stress singularities via semi-analytic separation of variables and natural compatibility with curvilinear coordinates \eqref{coordinate}, avoiding mesh generation difficulties for irregular geometries.

Denote $\mathscr{L} = (a_{11}^1, \cdots,a_{11}^m,a_{22}^1, \cdots,a_{22}^m,a_{33}^1, \cdots,a_{33}^m,a_{12}^1, \cdots,a_{12}^m,a_{13}^1, \cdots,a_{13}^m,$
$a_{23}^1, \cdots,a_{23}^m)^T,$ 
$
\mathcal{G} = ( \frac {\partial J}{\partial a_{11}^1}, \cdots, \frac {\partial J}{\partial a_{11}^m}, \frac {\partial J}{\partial a_{22}^1}, \cdots, \frac {\partial J}{\partial a_{22}^m},\frac {\partial J}{\partial a_{33}^1}, \cdots, \frac {\partial J}{\partial a_{33}^m},\frac {\partial J}{\partial a_{12}^1}, \cdots, \frac {\partial J}{\partial a_{12}^m},$
$\frac {\partial J}{\partial a_{13}^1}, \cdots, \frac {\partial J}{\partial a_{13}^m}, \frac {\partial J}{\partial a_{23}^1}, \cdots, \frac {\partial J}{\partial a_{23}^m})^T.
$
At the $k$-th iteration, we assign $\mathscr{L}^k$ as the current value of $\mathscr{L}$,   $a_h^k $ be the coefficients \eqref{discrete Lame coefficients} associated with $\mathscr{L}^k$, $\mathcal{G}^k$ represents $\mathcal{G}$ associated with the coefficients $a_h^k$. The Adam algorithm using the semi-analytical forward solver for heterogeneous and anisotropic elastic parameter estimation is similar to Algorithm 1 in \cite{zhu2023}.

\section{The direct method of lines as the forward solver}
\label{sec:GDMOL}

\subsection{The curvilinear coordinate}
\begin{figure}[H]
\centering
\includegraphics[scale=0.3]{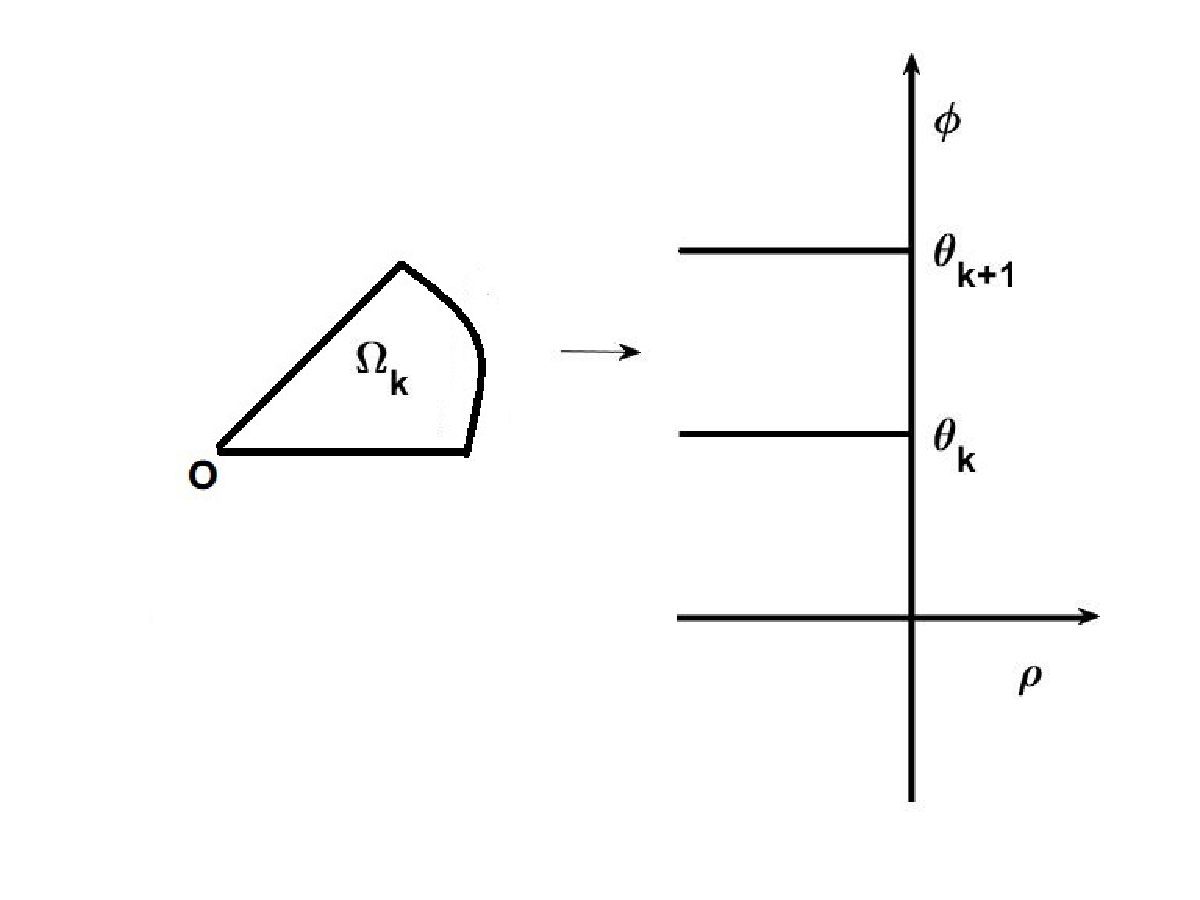}\vspace{-10mm}
\caption{Curvilinear coordinate remapping}
\label{fig:changedomain}
\end{figure}

The numerical solution of the composite elasticity problem \eqref{composite-eq} presents two fundamental challenges when using traditional Cartesian-based methods: (1)the star-shaped boundary cannot be efficiently discretized using regular grids, (2)piecewise-constant anisotropic coefficients with discontinuities across radial interfaces. We therefore introduce the coordinate system where
\begin{equation}
(x, y) = (e^{\rho} \tilde{r}(\phi) \cos(\phi), e^{\rho}\tilde{r}(\phi) \sin(\phi)),\quad -\infty < \rho \leq 0, \quad 0 \leq \phi < 2\pi.
\label{coordinate}
\end{equation}
Under mapping \eqref{coordinate}, each material subdomain $\Omega_k$ is transformed to a regular semi-infinite strip; see Figure \ref{fig:changedomain}. Thus more general irregular star-shaped domains can be solved.

Subsequently, transform \eqref{composite-eq} into its equivalent expression within the curvilinear coordinate system. This derivation begins with
$$
e^{\rho}\tilde{r}^2(\phi)\left(
\begin{matrix}
\frac{\partial \rho}{\partial x} & \frac{\partial \rho}{\partial y}\\
\frac{\partial \phi}{\partial x} & \frac{\partial \phi}{\partial y}
\end{matrix}\right)=\left(
\begin{matrix}
\tilde{r}^{\prime}(\phi) \sin(\phi)+\tilde{r}(\phi)\cos(\phi) & -\tilde{r}^{\prime}(\phi) \cos(\phi)+\tilde{r}(\phi)\sin(\phi)\\
-\tilde{r}(\phi) \sin(\phi) & \tilde{r}(\phi) \cos(\phi)
\end{matrix}\right).
$$
Then the stress tensor become
\begin{equation} 
\begin{aligned}
\sigma^k_{11} &=a_{11}^k\frac{\partial u_1^k}{\partial x}+a_{12}^k\frac{\partial u_2^k}{\partial y}+a_{13}^k(\frac{\partial u_1^k}{\partial y}+\frac{\partial u_2^k}{\partial x})\\
&=\frac1{e^{\rho}\tilde{r}^2}\bigg(a_{11}^k(\frac{\partial u_1^k}{\partial\rho}(\tilde{r}'\sin(\phi)+\tilde{r}\cos(\phi))+\frac{\partial u_1^k}{\partial\phi}(-\tilde{r}\sin{\phi}))+a_{12}^k(\frac{\partial u_2^k}{\partial\rho}\\
&(-\tilde{r}'\cos(\phi)+\tilde{r}\sin(\phi))+\frac{\partial u_2^k}{\partial\phi}(\tilde{r}\cos{\phi}))+a_{13}^k(\frac{\partial u_1^k}{\partial\rho}(-\tilde{r}'\cos(\phi)\\
&+\tilde{r}\sin(\phi))+\frac{\partial u_1^k}{\partial\phi}(\tilde{r}\cos{\phi})
+\frac{\partial u_2^k}{\partial\rho}(\tilde{r}'\sin(\phi)+\tilde{r}\cos(\phi))+\frac{\partial u_2^k}{\partial\phi}(-\tilde{r}\sin{\phi})\bigg),
\end{aligned} 
\end{equation}
\begin{equation} 
\begin{aligned}
\sigma^k_{22} &=a_{12}^k\frac{\partial u_1^k}{\partial x}+a_{22}^k\frac{\partial u_2^k}{\partial y}+a_{23}^k(\frac{\partial u_1^k}{\partial y}+\frac{\partial u_2^k}{\partial x})\\
&=\frac1{e^{\rho}\tilde{r}^2}\bigg(a_{12}^k(\frac{\partial u_1^k}{\partial\rho}(\tilde{r}'\sin(\phi)+\tilde{r}\cos(\phi))+\frac{\partial u_1^k}{\partial\phi}(-\tilde{r}\sin{\phi}))+a_{22}^k(\frac{\partial u_2^k}{\partial\rho}\\
&(-\tilde{r}'\cos(\phi)+\tilde{r}\sin(\phi))+\frac{\partial u_2^k}{\partial\phi}(\tilde{r}\cos{\phi}))+a_{23}^k(\frac{\partial u_1^k}{\partial\rho}(-\tilde{r}'\cos(\phi)\\
&+\tilde{r}\sin(\phi))+\frac{\partial u_1^k}{\partial\phi}(\tilde{r}\cos{\phi})+\frac{\partial u_2^k}{\partial\rho}(\tilde{r}'\sin(\phi)+\tilde{r}\cos(\phi))+\frac{\partial u_2^k}{\partial\phi}(-\tilde{r}\sin{\phi})\bigg),
\end{aligned}
\end{equation}
\begin{equation} 
\begin{aligned}
\sigma^k_{12} &=\sigma^k_{21}=a_{13}^k\frac{\partial u_1^k}{\partial x}+a_{23}^k\frac{\partial u_2^k}{\partial y}+a_{33}^k(\frac{\partial u_1^k}{\partial y}+\frac{\partial u_2^k}{\partial x})\\
&=\frac1{e^{\rho}\tilde{r}^2}\bigg(a_{13}^k(\frac{\partial u_1^k}{\partial\rho}(\tilde{r}'\sin(\phi)+\tilde{r}\cos(\phi))+\frac{\partial u_1^k}{\partial\phi}(-\tilde{r}\sin{\phi}))+a_{23}^k(\frac{\partial u_2^k}{\partial\rho}\\
&(-\tilde{r}'\cos(\phi)+\tilde{r}\sin(\phi))+\frac{\partial u_2^k}{\partial\phi}(\tilde{r}\cos{\phi}))+a_{33}^k(\frac{\partial u_1^k}{\partial\rho}(-\tilde{r}'\cos(\phi)\\
&+\tilde{r}\sin(\phi))+\frac{\partial u_1^k}{\partial\phi}(\tilde{r}\cos{\phi})+\frac{\partial u_2^k}{\partial\rho}(\tilde{r}'\sin(\phi)+\tilde{r}\cos(\phi))+\frac{\partial u_2^k}{\partial\phi}(-\tilde{r}\sin{\phi})\bigg).
\end{aligned}
\end{equation}
In $\Omega_k$, the original problem \eqref{composite-eq} can now be reformulated as :
\begin{equation} \label{discontinuous-composite-eq}
\begin{aligned}
&((\tilde{r}'sin(\phi)+\tilde{r}cos(\phi))\frac{\partial}{\partial \rho}-\tilde{r}sin(\phi)\frac{\partial}{\partial \phi})\sigma_{11}^k\\
&+((-\tilde{r}'cos(\phi)+\tilde{r}sin(\phi))\frac{\partial}{\partial \rho}+\tilde{r}cos(\phi)\frac{\partial}{\partial \phi})\sigma_{12}^k  = 0,\\
&((\tilde{r}'sin(\phi)+\tilde{r}cos(\phi))\frac{\partial}{\partial \rho}-\tilde{r}sin(\phi)\frac{\partial}{\partial \phi})\sigma_{12}^k\\
&+((-\tilde{r}'cos(\phi)+\tilde{r}sin(\phi))\frac{\partial}{\partial \rho}+\tilde{r}cos(\phi)\frac{\partial}{\partial \phi})\sigma_{22}^k  = 0,\\
 &u^k(0,\phi)=f^k(\phi), \quad u^{k - 1}(\rho,\theta_k^-) = u^{k}(\rho,\theta_k^+),\\
&\left(
\begin{matrix}
\sigma^{k - 1}_{11} \sin(\theta_k) - \sigma^{k-1}_{12} \cos(\theta_k) \\
\sigma^{k - 1}_{21} \sin(\theta_k) - \sigma^{k-1}_{22} \cos(\theta_k)
\end{matrix}\right)  = \left(
\begin{matrix}
\sigma^{k}_{11} \sin(\theta_k) - \sigma^{k}_{12} \cos(\theta_k) \\
\sigma^{k}_{21} \sin(\theta_k) - \sigma^{k}_{22} \cos(\theta_k)
\end{matrix}\right),\\
&u \text{ is bounded,} \text{ when } \rho\rightarrow -\infty. 
\end{aligned}
\end{equation}

 Denote $H^1_p(0, 2\pi)$ as the periodic Sobolev space, $V = H^1_p(0, 2\pi) \times H^1_p(0, 2\pi)$ and $W = \{w(\rho,\phi)|w, \partial_{\rho} w, \partial^2_{\rho\rho} w \in V , \forall\rho < 0 \}$.  Problem \eqref{discontinuous-composite-eq} can be equivalently expressed as
\begin{equation} \label{v-d-problem}
\begin{aligned}
& \text{determine } u\in W, \text{ s.t. }\forall v\in V,-\infty <\rho<0, \\
& \frac{d^2A_2(u,v)}{d\rho^2}+\frac{dA_1(u,v)}{d\rho}+A_0(u,v)=0,\\
& u(-\infty,\phi) <\infty, \quad u(0,\phi)=f(\phi),
\end{aligned}
\end{equation}
with
\begin{equation}
A_0(u,v) = -\sum_{k = 1}^K \int_{\theta_k}^{\theta_{k+1}} \left( \frac{\partial u^k (\rho, \phi)}{\partial \phi} \right)^T \Psi^k_0\frac{d v^k}{d\phi} d\phi,
\end{equation}
\begin{equation}
A_1(u,v) = \sum_{k = 1}^K \int_{\theta_k}^{\theta_{k+1}} \left( \left( \frac{\partial u^k (\rho, \phi) }{\partial \phi} \right)^T \Psi^k_1 v^k - u^k(\rho,\phi)^T (\Psi_1^k)^T \frac{d v^k}{d\phi} \right) d\phi,
\end{equation}
\begin{equation}
A_2(u,v) = \sum_{k = 1}^K \int_{\theta_k}^{\theta_{k+1}} u^k(\rho,\phi)^T \Psi_2^k v^k d\phi,
\end{equation}
where $v^k = v|_{[\theta_k,\theta_{k+1}]}$ and
$$
\begin{aligned}
\Psi_2^k = \frac1{\tilde{r}^2}\left(
\begin{matrix}
a_{11}^k(\tilde{r}'\sin(\phi)+\tilde{r}\cos(\phi))^2 & a_{13}^k(\tilde{r}'\sin(\phi)+\tilde{r}\cos(\phi))^2\\
+2a_{13}^k(\tilde{r}'\sin(\phi)+\tilde{r}\cos(\phi)) & +(a_{33}^k+a_{12}^k)(\tilde{r}'\sin(\phi)+\tilde{r}\cos(\phi))\\
(-\tilde{r}'\cos(\phi)+\tilde{r}\sin(\phi)) & (-\tilde{r}'\cos(\phi)+\tilde{r}\sin(\phi))\\
+a_{33}^k(-\tilde{r}'\cos(\phi)+\tilde{r}\sin(\phi))^2 &+a_{23}^k(-\tilde{r}'\cos(\phi)+\tilde{r}\sin(\phi))^2\\
\\
a_{13}^k(\tilde{r}'\sin(\phi)+\tilde{r}\cos(\phi))^2 &  a_{33}^k(\tilde{r}'\sin(\phi)+\tilde{r}\cos(\phi))^2\\
+(a_{33}^k+a_{12}^k)(\tilde{r}'\sin(\phi)+\tilde{r}\cos(\phi))& +2a_{23}^k(\tilde{r}'\sin(\phi)+\tilde{r}\cos(\phi)) \\
(-\tilde{r}'\cos(\phi)+\tilde{r}\sin(\phi))& (-\tilde{r}'\cos(\phi)+\tilde{r}\sin(\phi))\\
+a_{23}^k(-\tilde{r}'\cos(\phi)+\tilde{r}\sin(\phi))^2&  +a_{22}^k(-\tilde{r}'\cos(\phi)+\tilde{r}\sin(\phi))^2
\end{matrix}\right), 
\end{aligned}
$$
$$
\Psi_1^k = \frac1{\tilde{r}}\left(
\begin{matrix}
-a_{11}^k(\tilde{r}'\sin(\phi)+\tilde{r}\cos(\phi))\sin(\phi) & -a_{13}^k(\tilde{r}'\sin(\phi)+\tilde{r}\cos(\phi))\sin(\phi) \\
+a_{13}^k(\tilde{r}'\sin(\phi)+\tilde{r}\cos(\phi))\cos(\phi)& +a_{12}^k(\tilde{r}'\sin(\phi)+\tilde{r}\cos(\phi))\cos(\phi)\\
-a_{13}^k(-\tilde{r}'\cos(\phi)+\tilde{r}\sin(\phi))\sin(\phi)&-a_{33}^k(-\tilde{r}'\cos(\phi)+\tilde{r}\sin(\phi))\sin(\phi)\\
+a_{33}^k(-\tilde{r}'\cos(\phi)+\tilde{r}\sin(\phi))\cos(\phi)&+a_{23}^k(-\tilde{r}'\cos(\phi)+\tilde{r}\sin(\phi))\cos(\phi)\\
\\
-a_{13}^k(\tilde{r}'\sin(\phi)+\tilde{r}\cos(\phi))\sin(\phi) & -a_{33}^k(\tilde{r}'\sin(\phi)+\tilde{r}\cos(\phi))\sin(\phi) \\
+a_{33}^k(\tilde{r}'\sin(\phi)+\tilde{r}\cos(\phi))\cos(\phi)& +a_{23}^k(\tilde{r}'\sin(\phi)+\tilde{r}\cos(\phi))\cos(\phi)\\
-a_{12}^k(-\tilde{r}'\cos(\phi)+\tilde{r}\sin(\phi))\sin(\phi)&-a_{23}^k(-\tilde{r}'\cos(\phi)+\tilde{r}\sin(\phi))\sin(\phi)\\
+a_{23}^k(-\tilde{r}'\cos(\phi)+\tilde{r}\sin(\phi))\cos(\phi)&+a_{22}^k(-\tilde{r}'\cos(\phi)+\tilde{r}\sin(\phi))\cos(\phi)
\end{matrix}\right), 
$$
$$
\Psi_0^k = \left(
\begin{matrix}
a_{11}^k\sin^2(\phi)+a_{33}^k\cos^2(\phi)& a_{13}^k\sin^2(\phi)+a_{23}^k\cos^2(\phi)\\
-2a_{13}^k\cos(\phi)\sin(\phi)& -(a_{33}^k+a_{12}^k)\cos(\phi)\sin(\phi) \\
\\
a_{13}^k\sin^2(\phi)+a_{23}^k\cos^2(\phi)& a_{33}^k\sin^2(\phi)+a_{22}^k\cos^2(\phi)\\
-(a_{33}^k+a_{12}^k)\cos(\phi)\sin(\phi)&-2a_{23}^k\cos(\phi)\sin(\phi)
\end{matrix}\right).
$$
\begin{lemma}\label{lem:A} The following properties hold for the three bounded bilinear forms $A_j$ on $V\times V$:\\
(1) $A_0$ and $A_2$ are symmetric, while $A_1$ is antisymmetric.\\
(2) There exists a constant $c>0$ such that
\[
-A_0(v,v)\geq c|| v^{\prime} ||_2^2,\quad A_2(v,v)\geq c || v ||_2^2, \quad \forall v \in V.
\]
where $|| \cdot ||_2$ is the $L^2$ norm and $v^{\prime}$ is the weak derivative of $v$. 
\end{lemma}

\subsection{Semi-analytical treatment}
To numerically solve  variational-differential problem \eqref{v-d-problem},
let $[0,2\pi]$ be partitioned into $M$ subintervals with nodes $0 = \phi_1 < \ldots < \phi_{M + 1} = 2\pi$, such that all material interface angles $\{\theta_k\}_{k = 1}^K$ are included in the nodal set and $h = \max_{1 \le j \le M} |\phi_{j+1}-\phi_j|$.  Choose $V_1^h\subset H^1_p(0, 2\pi)$ to be $P_1$-conforming linear finite element space, $V^h = V_1^h \times V_1^h$ and 
$W^h = \left\{ w(\rho,\phi)|w, \frac{\partial w}{\partial\rho} , \frac{\partial^2 w}{\partial\rho^2}  \in V^h, \forall\rho < 0 \right\}.$
Next, we present the semi-discrete approximation of \eqref{v-d-problem} as
\begin{equation} \label{n-v-d-problem}
\begin{aligned}
& \text{determine } u^h\in W^h, \text{ s.t. }\forall v^h\in V^h,-\infty <\rho<0, \\
& \frac{d^2A_2(u^h,v^h)}{d\rho^2}+\frac{dA_1(u^h,v^h)}{d\rho}+A_0(u^h,v^h)=0,  \\
&  u^h(-\infty,\phi) < \infty, \quad u^h(0,\phi)=f^h(\phi),
\end{aligned}
\end{equation}
with $f^h(\phi) \in V^h$ s.t. $f^h(\phi_j) = f(\phi_j)$ for $j = 1, \ldots, M$. 
Denote the basis functions of $V_1^h$ be $\{BF_j(\phi)\}^M_{j=1}$ s.t. $BF_j(\phi_k)=\delta_{kj}, 1\leq k,j\leq M$ and
\[
BF(\phi)=\left(
\begin{matrix}
BF_1(\phi) & BF_2(\phi) &... & BF_M(\phi) & 0 & 0&...& 0\\
0 & 0&...&0  & BF_1(\phi) & BF_2(\phi) &... & BF_M(\phi)
\end{matrix}\right)^T.
\]
Let the semi-analytical solution be
\begin{equation} \label{semi approximation of the solution}
u^h(\rho,\phi)=BF^T(\phi)U(\rho).
\end{equation}
If
$
F=(f_1(\phi_1),...,f_1(\phi_M),f_2(\phi_1),...,f_2(\phi_M))^T,
$
then
$
f^h(\phi)=BF^T(\phi)F.
$
Therefore, \eqref{n-v-d-problem} reduces to the following second-order ODEs
\begin{equation} \label{eq:odes}
\begin{aligned}
& B_2U''(\rho)+B_1U'(\rho)+B_0U(\rho)=0,\quad \forall\rho\in(-\infty,0), \\
& U(-\infty) < \infty,\quad U(0)=F, 
\end{aligned}
\end{equation}
with
\begin{gather}
B_0 = -\sum_{k=1}^K \int_{\theta_k}^{\theta_{k+1}} BF'(\phi)\Psi_0^k BF'(\phi)^T d\phi,\\
B_1 = \sum_{k=1}^K \int_{\theta_k}^{\theta_{k+1}} \left( BF(\phi) \Psi_1^k BF'(\phi)^T - BF'(\phi) (\Psi_1^k)^T BF(\phi)^T \right) d\phi,\\
B_2 = \sum_{k=1}^K \int_{\theta_k}^{\theta_{k+1}} BF(\phi) \Psi_2^k BF(\phi)^T d\phi,
\end{gather}
where $BF'(\phi)$ be the derivative of $BF(\phi)$.  Then 
\begin{lemma} \label{lem: properties of B}
The constant matrix $B_2$ exhibits symmetric positive definiteness, whereas constant matrix $B_1$ demonstrates antisymmetry and constant matrix $B_0$ possesses symmetric negative semi-definite characteristics.
\end{lemma}
Which derived from  Lemma \ref{lem:A}.

We employ a direct method to solve the system of ODEs \eqref{eq:odes}. Let
\begin{equation} \label{form of U}
U(\rho)=e^{\rho\gamma}\xi,
\end{equation}
with $\gamma$ being an unknown constant and $\xi\in C^{2M}$ an unknown vector function. By inserting the solution form \eqref{form of U} into the ODEs \eqref{eq:odes}, we derive the quadratic eigenvalue problem
\begin{equation}
[\gamma^2B_2+\gamma B_1+B_0]\xi=0.
\label{eq:eig-pro}
\end{equation}
By introducing the auxiliary variable $\zeta=\gamma\xi$, an  linear generalized eigenvalue problem equivalent  to \eqref{eq:eig-pro} yields
\begin{equation}
\left(
\begin{matrix}
0 & I_{2M}\\
-B_0 & -B_1
\end{matrix}\right)\left(
\begin{matrix}
\xi\\
\zeta
\end{matrix}\right)=\gamma\left(
\begin{matrix}
I_{2M} & 0\\
0 & B_2
\end{matrix}\right)\left(
\begin{matrix}
\xi\\
\zeta
\end{matrix}\right),
\label{eq:general-eig}
\end{equation}
where $I_{2M}$ is the identity matrix. Combine the general theory of quadratic eigenvalue problems developed in \cite{quadratic-eig} with Lemma \ref{lem: properties of B}, we derive
\begin{lemma}
The generalized eigenvalue problem (\ref{eq:general-eig}) possesses exactly $2M$ eigenvalues with $Re(\gamma)\geq 0$ and $2M$ eigenvalues with $Re(\gamma)\leq 0$, counting multiplicities.
\end{lemma}

For the generalized eigenvalue problem \eqref{eq:general-eig}, numerical computation gives $2M$ eigenvalues $\{\gamma_j^h\}^{2M}_{j=1}$ with $Re(\gamma_j^h)\geq 0$ and corresponding eigenvectors $(\xi_j^T, \zeta_j^T)^T$, $j = 1, \ldots, 2M$.
The problem admits at least two zero eigenvalues, which we denote as $\gamma_1^h = \gamma_2^h = 0$. Let $\{\gamma_j^h\}_{j=1}^{2m}$ be the real eigenvalues (ordered ascendingly) and $\{\gamma_j^h\}_{j=2m+1}^{2M}$ be the complex eigenvalues with $Im(\gamma_j^h)\neq 0$, satisfying $\gamma_{2j}^h = \overline \gamma_{2j-1}^h,\, m+1 \leq j \leq M$. Hence
\begin{equation}\label{eq:U}
U(\rho) = \sum_{j=1}^{2m} \alpha_j e^{\rho \gamma_j^h} \xi_j + \sum_{j = 1+m}^{M} \left( \alpha_{2j-1} \textbf{Re} (e^{\rho\gamma_{2j}^h} \xi_{2j}) + \alpha_{2j} \textbf{Im} (e^{\rho\gamma_{2j}^h} \xi_{2j}) \right).
\end{equation}
 Since $U(0)=F$, we obtain
\begin{equation} \label{eq:F}
F = \sum_{j=1}^{2m} \alpha_j \xi_j + \sum_{j=1+m}^{M} \left( \alpha_{2j-1} \textbf{Re}(\xi_{2j}) + \alpha_{2j} \textbf{Im}(\xi_{2j}) \right).
\end{equation}
Let
$$
\begin{aligned}
EI(\rho) = [ & e^{\rho\gamma_1^h}\xi_1 , \ldots, e^{\rho\gamma_{2m}^h} \xi_{2m}, \textbf{Re}(e^{\rho\gamma_{2m+2}^h}\xi_{2m+2}), \textbf{Im}(e^{\rho\gamma_{2m+2}}\xi_{2m+2}), \ldots, \\
& \textbf{Re}(e^{\rho\gamma_{2M}^h}\xi_{2M}), \textbf{Im}(e^{\rho\gamma_{2M}^h}\xi_{2M})],
\end{aligned}
$$
$$
EI(0)=[\xi_1,...,\xi_{2m}, \textbf{Re}(\xi_{2m+2}), \textbf{Im}(\xi_{2m+2}), \ldots, \textbf{Re}(\xi_{2M}), \textbf{Im}(\xi_{2M})],
$$
$$
\alpha = [\alpha_1, \alpha_2, \ldots, \alpha_{2M}]^T.
$$
Equation (\ref{eq:F}) yields
\begin{equation}\label{eq:BB}
\alpha = EI(0)^{-1}F.
\end{equation}
Inserting (\ref{eq:BB}) into (\ref{eq:U}) gives
\begin{equation}
U(\rho) = EI(\rho) EI(0)^{-1} F.
\end{equation}
The variational problem \eqref{v-d-problem} admits a semi-discrete solution through the approximation scheme \eqref{semi approximation of the solution}, giving
\begin{equation}\label{semi-dis-solution}
u^h(\rho,\phi) = BF(\phi)^T EI(\rho) EI(0)^{-1} F.
\end{equation}

For clarity in demonstrating how our method captures solution singularities, we first consider the case where all eigenvalues $\{\gamma_j^h\}^{2M}_{j=1}$ of the quadratic eigenvalue problem \eqref{eq:eig-pro} are real. The extension to cases with complex eigenvalues follows analogous arguments. From the semi-discrete formulation \eqref{semi-dis-solution}, our numerical solution is
$
u^h = \sum_{j = 1}^{2M} \alpha_j r^{\gamma^h_j} (\Phi^h_j(\phi), \Psi^h_j(\phi))^T
$
in the polar coordinate, with
$$
\Phi^h_j(\phi) = \left( \tilde{r}(\phi) \right)^{-\gamma^h_j} \sum_{i = 1}^M BF_k(\phi) (\xi_j)_i, \quad \Psi^h_j(\phi) = \left( \tilde{r}(\phi) \right)^{-\gamma^h_j} \sum_{i = 1}^M BF_k(\phi) (\xi_j)_{M + i}.
$$
Our numerical experiments show fast convergence of $\gamma^h_j$ to the elliptic operator's true eigenvalue $b_j$, indicating that the method naturally handles singularities—refer to Section \ref{subsec: numerical examples of forward problems}. Let $u$ be the solution to \eqref{composite-eq}. The error estimate can be derived using arguments analogous to those in Theorem 4.1 of \cite{com-material}.
\begin{thm}\label{thm-error}
If $V_1^h$ consists of linear elements, then there exists $C > 0$ independent of the mesh size $h$ for which
\begin{equation} \label{eq: error estimate}
\|u-u^h\|^2_*\leq Ch^2\sum_{j=1}^2\iint_{\Omega} ( |\nabla u_j|^2+\sum_{\beta_1+\beta_2=2}(x^2+y^2)|\frac{\partial^{2} u_j}{\partial x^{\beta_1} \partial y^{\beta_2}}|^2 ) dxdy.
\end{equation}
where
\begin{gather}
\|v\|_* = E(v,v)^{1/2} + |\psi_1(v)| + |\psi_2(v)| + |\psi_3(v)|, \notag\\
\psi_1(v) = \int_{\Gamma}v_1ds,\quad
\psi_2(v) = \int_{\Gamma}v_2ds,\quad
\psi_3(v) = \int_{\Gamma}(v_1y-v_2x)ds,\notag
\end{gather}

$$
\begin{aligned}
E(w,v) & = \iint_{\Omega} \bigg( a_{11}\frac{\partial w_1}{\partial x} \frac{\partial v_1}{\partial x}+a_{22}\frac{\partial w_2}{\partial y}\frac{\partial v_2}{\partial y}+a_{33}(\frac{\partial w_1}{\partial y}+\frac{\partial w_2}{\partial x})(\frac{\partial v_1}{\partial y}+\frac{\partial v_2}{\partial x}) \\
&+a_{12}\frac{\partial w_1}{\partial x}\frac{\partial v_2}{\partial y}+a_{21}\frac{\partial v_1}{\partial x}\frac{\partial w_2}{\partial y}+ a_{13}\frac{\partial w_1}{\partial x}(\frac{\partial v_2}{\partial x}+\frac{\partial v_1}{\partial y})+a_{31}\frac{\partial v_1}{\partial x}(\frac{\partial w_2}{\partial x}+\frac{\partial w_1}{\partial y})\\
&+a_{23}\frac{\partial w_2}{\partial y}(\frac{\partial v_2}{\partial x}+\frac{\partial v_1}{\partial y}) +a_{32}\frac{\partial v_2}{\partial y}(\frac{\partial w_2}{\partial x}+\frac{\partial w_1}{\partial y}) \bigg) dxdy.
\end{aligned}
$$
\end{thm}

Observe that while the solution $u$ belongs to $H^1(\Omega)$, it does not possess  $H^2(\Omega)$-regularity. As a result, $|\nabla u|^2$ is finite, but $|D^{\beta} u|^2$ diverges for all second-order derivatives. However, the weighted term $r^2 |D^{\beta} u|^2$ is integrable because $u$ develops singularities solely along the radial direction at the origin, with asymptotic structure $\mathcal{O}(r^{b_j})$ for $0 < \textbf{Re}(b_j) < 1$. Theorem \ref{thm-error} thus guarantees that the right-hand side of \eqref{eq: error estimate} is finite, yielding first-order accuracy with respect to the $|| \cdot ||_{*}$-norm, even when $u$ lacks $H^2(\Omega)$-regularity. Combined with the symmetric positive definiteness of $a$, we derive the 
$H^1$-error control $|| u - u^h ||_{H^1} \le C || u - u^h ||_{*}$ with $C$ independent of discretization parameters. Consequently, when the tensor $a$ is symmetric positive definite, linear partial discrete approximation achieves first-order $H^1$-norm accuracy, as proved in Theorem \ref{thm-error}. This result naturally extends to quadratic elements achieve second-order convergence, with numerical confirmation provided in Example \ref{ex2}.

\section{Numerical examples} \label{sec:example}

\subsection{The forward problems} \label{subsec: numerical examples of forward problems}
In Example \ref{ex2}, we conduct a direct comparison between linear and quadratic elements for semi-discretization, employing a known exact solution for verification. And in Examples \ref{ex1}, we generate numerical solutions and approximate eigenvalues using linear elements, while establishing reference values through high-resolution computations with quadratic elements. We measure errors using the relative $H^1$-norm $\frac{||u-u^h||_{1}}{||u||_{1}}$.

\begin{figure}[H]
\centering
\includegraphics[scale=0.4]{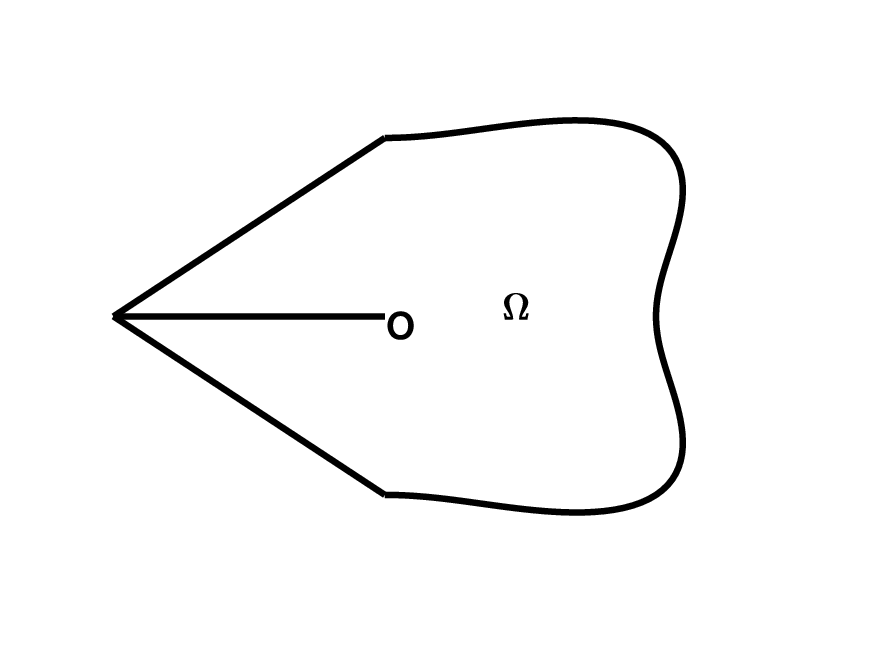}
\caption{the irregular stellar region of Example \ref{ex2}}
\label{fig:domianExp2}
\end{figure}

\begin{example}\label{ex2} We parameterize the boundary as
$$
\tilde{r}(\phi) =
\begin{cases}
- \frac{1}{\sin(\phi)+\cos(\phi)}, & \quad  -\pi \le \phi \le -\pi/2,\\
\frac1{ \sin^4(\phi)+\cos^4(\phi)}, & \quad  -\pi/2 < \phi \le \pi/2,\\
- \frac{1}{-\sin(\phi)+\cos(\phi)}, & \quad \pi/2 < \phi \le \pi.
\end{cases}
$$
In Figure \ref{fig:domianExp2}, all points $(r,\phi)$ with $r\in(0, \tilde{r}(\phi))$ and $\phi\in(-\pi,\pi)$ comprise $\Omega$, boundary points with $\phi = -\pi, \pi$ compose $\Gamma_N$ and other boundary points constitute $\Gamma_D$.
\begin{equation} \label{equatioin of ex2}
\begin{aligned}
-\nabla\cdot \sigma & = 0,\quad u|_{\Gamma_D} & = f, \quad\sigma \cdot n|_{\Gamma_N} & = 0,
\end{aligned}
\end{equation}
with anisotropic coefficient tensor $a=\left(
\begin{matrix}
4.8 & 2.1 & 0.3\\
2.1 & 3.3 & 0.2\\
0.3 & 0.2 & 1.2
\end{matrix}\right)$.
 Let $b=a^{-1}$ and $\mu_1,\mu_2$ are the different complex roots of the equation
$$
b_{11}\mu^4-2b_{31}\mu^3+(2b_{21}+b_{33})\mu^2-2b_{32}\mu+b_{22}=0,
$$
$$
p_j=b_{11}\mu_j^2+b_{21}-b_{31}\mu_j,\quad q_j=b_{21}\mu_j+\frac{b_{22}}{\mu_j}-b_{32},\quad j=1,2.
$$
Then the exact solution\cite{Steigemann2015} of \eqref{equatioin of ex2} in the polar coordinate is
$$
\sqrt{\frac{2r}{\pi}}Re
\left(\frac1{\mu_1-\mu_2}(\mu_1\sqrt{\cos(\phi)+\mu_2\sin(\phi)}\left(
\begin{aligned}
p_2\\
q_2
\end{aligned}
\right)-\mu_2\sqrt{\cos(\phi)+\mu_1\sin(\phi)}\left(
\begin{aligned}
p_1\\
q_1
\end{aligned}
\right) )\right)
$$
 and stress singularities occur at the origin.
\end{example}

Think a uniform discretization of $[-\pi, \pi]$ with $h = \frac{2\pi}{M}$ yielding the grid points $-\pi = \phi_1 < \phi_2 < \cdots < \phi_{M + 1} = \pi$. Table \ref{tab:unit} presents the errors in approximating the smallest positive eigenvalue $\gamma_3^h$ with first-order and second-order finite element implementations for varying $M$. The results demonstrate that 
linear discrete approximation convergence with order two, whereas quadratic discrete approximation convergence with order four. This confirms our method's capability to naturally handle solution singularities. Table \ref{tab:errors-unit} presents the relative inaccuracy for computational approximations $u^h$ with first-order and second-order finite element implementations. The results demonstrate that first-order convergence in the $H^1$ relative error for linear elements and second-order convergence in the $H^1$ relative error for quadratic elements. The observed optimal convergence rates persist despite stress singularities, confirming our method's effectiveness for anisotropic elasticity problems. Moreover, higher-order elements yield accelerated convergence, as expected from approximation theory.

\begin{table}
\centering
\caption{In Example \ref{ex2} numerical inaccuracy for the smallest positive eigenvalue $\gamma_3^h$ }
\label{tab:unit}
\footnotesize
\begin{tabular}{@{}lllll}
\br
M & linear elements & convergence order & quadratic elements & convergence order \\
\mr
16 & 1.911e-3 &  & 4.883e-4 &  \\
32 & 3.150e-4 & 2.601 & 2.712e-5 & 4.170 \\
64 & 7.149e-5 & 2.139 & 1.750e-6 & 3.954 \\
128 & 1.741e-5 & 2.038 & 1.104e-7 & 3.987 \\
256 & 4.325e-6 & 2.010 & 6.913e-9 & 3.997 \\
\br
\end{tabular}
\end{table}

\begin{table}
\centering
\caption{In Example \ref{ex2} $H^1$ relative inaccuracy for $u^h$}
\label{tab:errors-unit}
\footnotesize
\begin{tabular}{@{}lllll}
\br
M & linear elements & convergence order & quadratic elements & convergence order\\
\mr
8 & 1.556e-1 &  & 7.591e-2 &  \\
16 & 9.527e-2 & 0.708 & 1.405e-2 & 2.434 \\
32 & 3.844e-2 & 1.309 & 3.042e-3 & 2.208 \\
64 & 1.616e-2 & 1.250 & 7.657e-4 & 1.990 \\
128 & 7.597e-3 & 1.089 & 1.922e-4 & 1.994 \\
\br
\end{tabular}
\end{table}

\begin{figure}[H]
\centering
\includegraphics[scale=0.4]{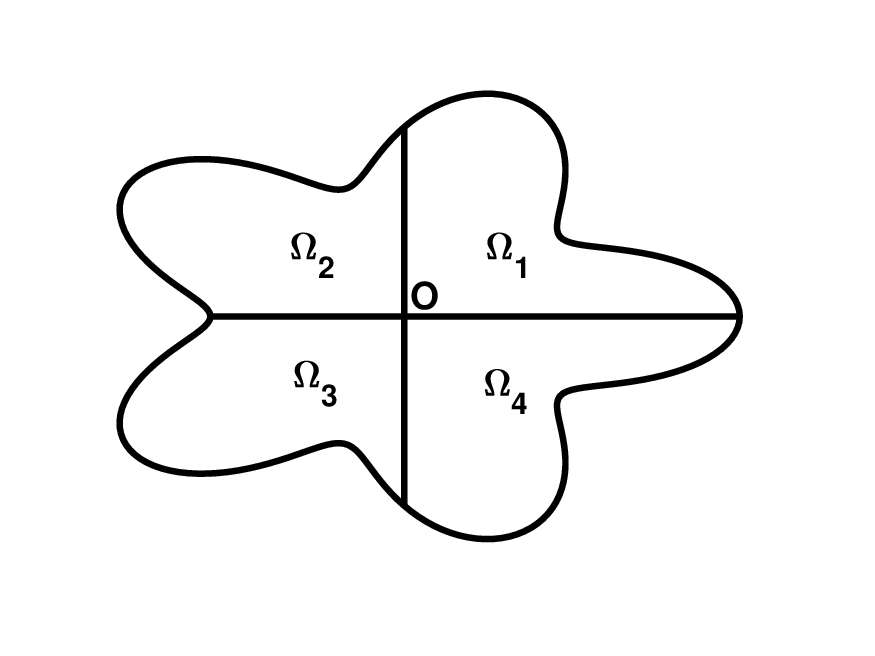}
\caption{Domain $\Omega$ in Example \ref{ex1}}
\label{fig:domainExp1}
\end{figure}

\begin{figure}[H]
\centering
\begin{subfigure}[b]{0.45\textwidth}
\includegraphics[width=\textwidth]{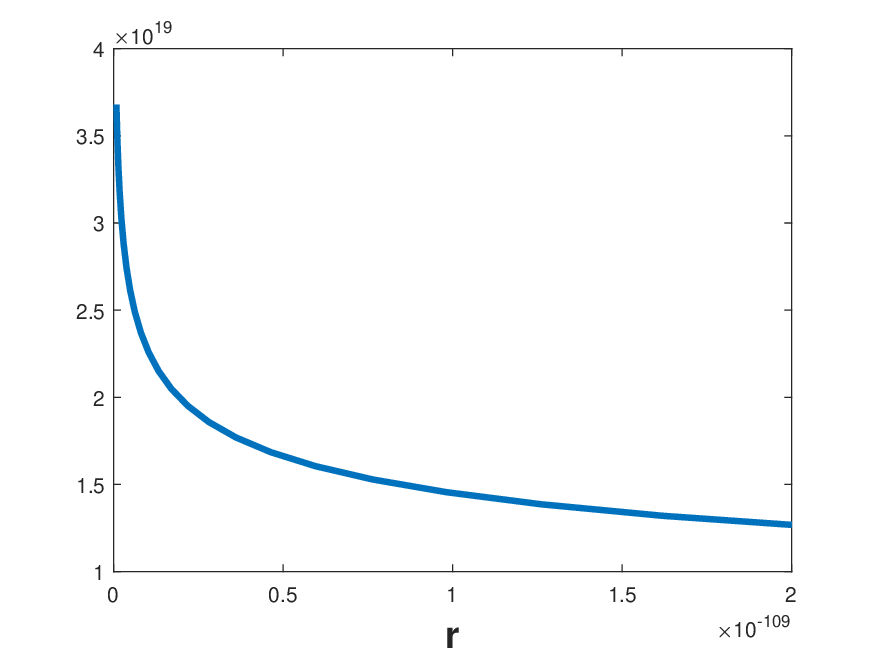}
\caption{$\frac{\partial u_1}{\partial r}$}
\end{subfigure}
\begin{subfigure}[b]{0.45\textwidth}
\includegraphics[width=\textwidth]{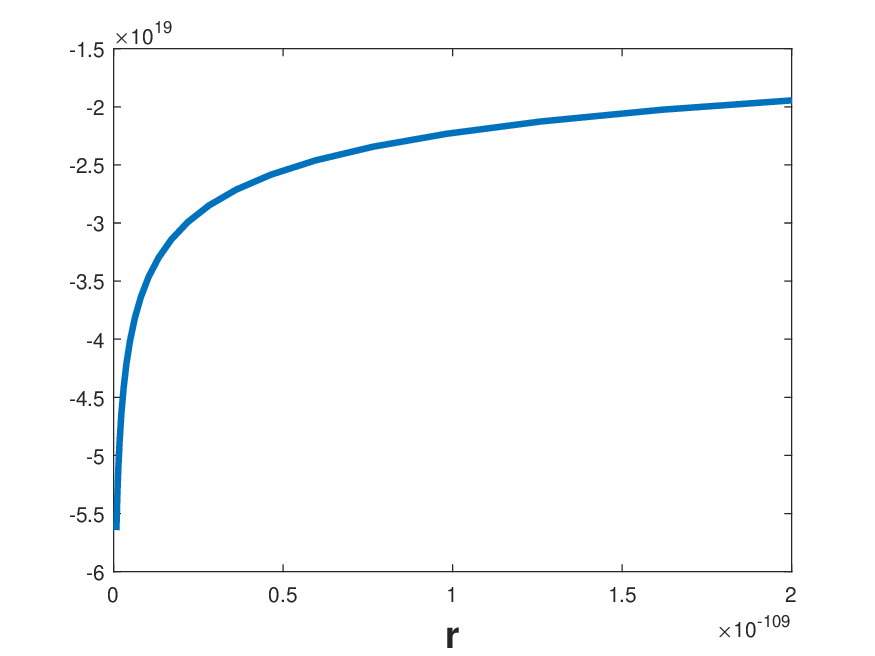}
\caption{$\frac{\partial u_2}{\partial r}$}
\end{subfigure}
\caption{$\frac{\partial u}{\partial r}(r, \phi_0)$ with $\phi_0 = \pi$ in Example \ref{ex1}}
\label{singularity of ex1}
\end{figure}

\begin{example}\label{ex1}
We parameterize the boundary as $\tilde{r}(\phi) = \sqrt{2 + \cos(5\phi)}$. All points $(r,\phi)$ with $r\in(0, \tilde{r}(\phi))$ and $\phi\in[0,2\pi)$ comprise $\Omega$, all points $(r,\phi)$ with $r\in(0, \tilde{r}(\phi))$ and $\phi\in(\theta_k,\theta_{k+1})$ comprise $\Omega_k$ where $\theta_1 = 0, \theta_2 = \pi/2, \theta_3 = \pi, \theta_4 = 3\pi/2, \theta_5 = 2\pi$; see Figure \ref{fig:domainExp1}. Consider
\begin{equation}
\begin{aligned}
-\nabla\cdot \sigma^k & = p, \quad \text{in } \Omega_k,  \\
u|_{\partial \Omega} & = f, \\
u^{k-1}|_{\theta_k^-} & = u^{k}|_{\theta_k^+},  \\
\sigma^{k-1}\cdot {n}_k|_{\theta_k^-} & = \sigma^{k}\cdot {n}_k|_{\theta_k^+}, 
\end{aligned}
\end{equation}
with $a_1=10a=10\left(
\begin{matrix}
4 & 2 & 0.2\\
2 & 3 & 0.1\\
0.2 & 0.1 & 1
\end{matrix}\right), a_2=5a, a_3=a, a_4=5a$, $p=(1,1)^T$ and $f=(1,1)^T$. As evidenced in Figure \ref{singularity of ex1}, the radial derivatives at $\phi_0 = \pi$ exhibit singular behavior near the origin: $\frac{\partial u_1}{\partial r}$ grows unboundedly to $\infty$ and $\frac{\partial u_2}{\partial r}$ diverges to  $-\infty$. This provides numerical verification of stress singularities in the solution at $r=0$.
\end{example}
The angular domain $[0,2\pi]$ is divided into $M$ equal sectors with $h = \frac{2\pi}{M}$. Table \ref{tab:multi-material} presents the approximation errors for the first nonzero eigenvalue $\gamma_3^h$, with reference value $\gamma_3 = 0.80480727808626$. The results demonstrate that second-order convergence for the eigenvalue approximation, while the $H^1$ relative inaccuracy for $u^h$ convergence with order one. The optimal eigenvalue convergence confirms our method's capability to naturally capture solution singularities, while the $H^1$ error convergence agrees with Theorem \ref{thm-error}.

\begin{table}
\centering
\caption{In Example \ref{ex1} inaccuracy of $\gamma_3^h$ and $H^1$ relative inaccuracy for $u^h$ }
\label{tab:multi-material}
\footnotesize
\begin{tabular}{@{}lllll}
\br
M & $| \gamma_3^h - \gamma_3 |$ & convergence order &$||u-u^h||_{1}/||u||_{1}$ & convergence order\\
\mr
8 & 2.880e-1 &  &1.156e-2 &\\
16 & 1.465e-1 & 0.975 &7.455e-3 &0.633 \\
32 & 2.216e-2 & 2.725 & 3.185e-3& 1.227\\
64 & 5.758e-3 & 1.944 & 1.448e-3& 1.137\\
128 & 1.455e-3 & 1.984 &6.991e-4 & 1.051\\
\br
\end{tabular}
\end{table}

Numerical experiments (Examples \ref{ex2}-\ref{ex1}) demonstrate that rapid convergence of discrete eigenvalues to their exact counterparts and optimal convergence rates for both linear and quadratic elements. These results confirm our method's effectiveness in naturally resolving stress singularities while maintaining accuracy for forward problems of anisotropic composites in star-shaped domains.

\subsection{The inverse problems} \label{subsec: numerical examples of inverse problems}
We employ  Adam Algorithm to solve the minimization problem \eqref{num-minpro} over the discrete parameter space $\Lambda_h$ with hyperparameter $\hat{\beta}_1 = 0.9$ and $\hat{\beta}_2 = 0.999$. Assume the regularization coefficient in problem \eqref{num-minpro} is $\eta = 1e-7$, the learning rate $\tau_k$ is diminishing, while the tolerance threshold is $tol =1e-5$. The anisotropic coefficient tensor $a$ is discretized using angular mesh size $h=2\pi/m$  and is represented as
$
a_h = \mathbbm{1}_{(-\infty, 0]}(\rho) \sum_{j=1}^{m} a_j \mathbbm{1}_{[(j-1)h, jh]}(\phi).
$
At each iteration $k$, we denote the parameter values by $\{a_j^k\}_{j=1}^m$. With mesh size $\tilde{h}=\pi/64$, the forward problem is solved via the semi-analytical method under linear discrete approximation. While reference solution $u[a^{\star}]$ is computed under quadratic finite elements via the semi-analytical method on a refined mesh. Measurement data is collected at discrete sensor locations $\Xi = \{ (x_j, y_j) \}_{j = 1}^{M_1} \subset \Omega$, where 
$M_1$ denotes the number of measurement point. The reconstruction accuracy is quantified through the coefficient error in $L^1$ norm.

\begin{figure}[H]
\centering
\includegraphics[scale=0.4]{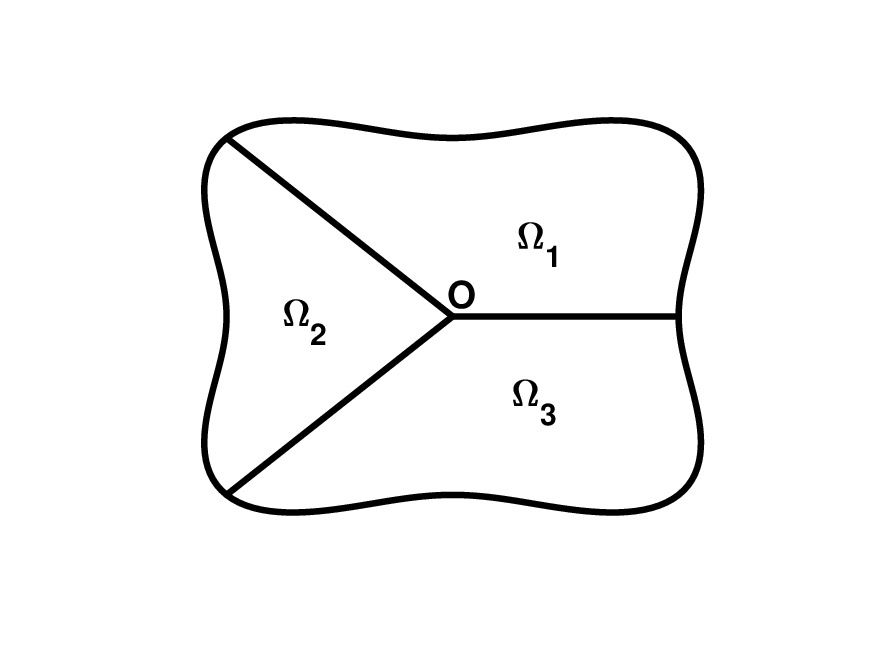}
\caption{Domain $\Omega$ in Example \ref{example1}} 
\label{fig:domainExample1}
\end{figure}

\begin{example}\label{example1}
Let $\tilde{r}(\phi) = (\cos^4(\phi) + \sin^4(\phi))^{-\frac 12}$, $\theta_1=0,\theta_2=3\pi/4$ and $\theta_3=5\pi/4$; see Figure \ref{fig:domainExample1}. 
Set $a_1^{\star}=\left(
\begin{matrix}
6 & 1 & 1\\
1 & 5 & 1\\
1 & 1 & 4
\end{matrix}\right)$,$a_2^{\star}=\left(
\begin{matrix}
5 & 1 & 1\\
1 & 4 & 1\\
1 & 1 & 3
\end{matrix}\right)$,$a_3^{\star}=\left(
\begin{matrix}
4 & 1 & 1\\
1 & 3 & 1\\
1 & 1 & 2
\end{matrix}\right)$ in distinct subdomains, Dirichlet boundary condition $f=(1, 1)^T$ and source term $p =(\frac{x}{\sqrt{x^2+y^2}},\frac{y}{\sqrt{x^2+y^2}})^T$.

\end{example}

\begin{figure}[H]
\minipage{0.45\textwidth}
\includegraphics[width=\linewidth]{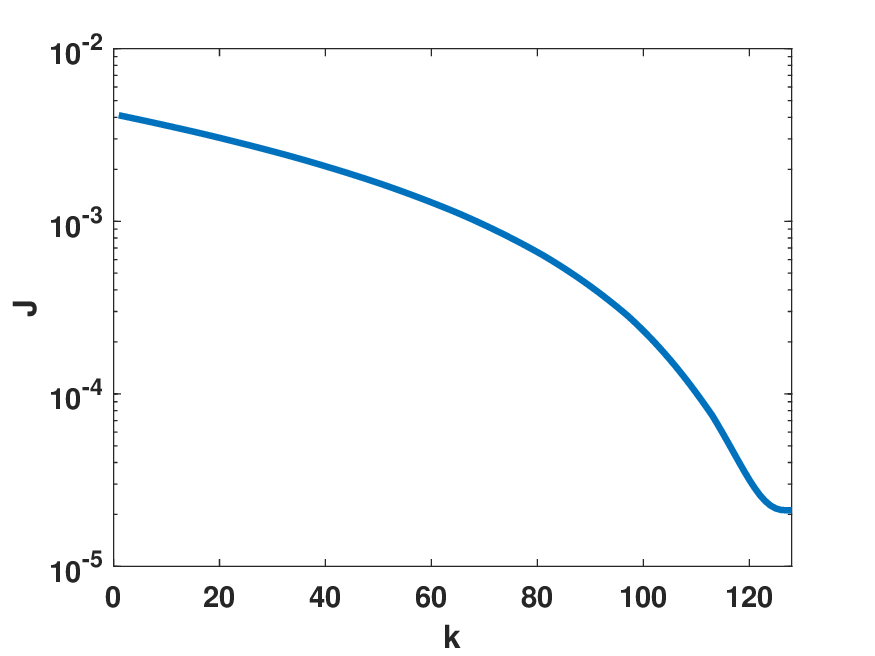}
\caption{Values of J in Example \ref{example1}} \label{fig::J_1}
\endminipage\hfill
\minipage{0.45\textwidth}
\includegraphics[width=\textwidth]{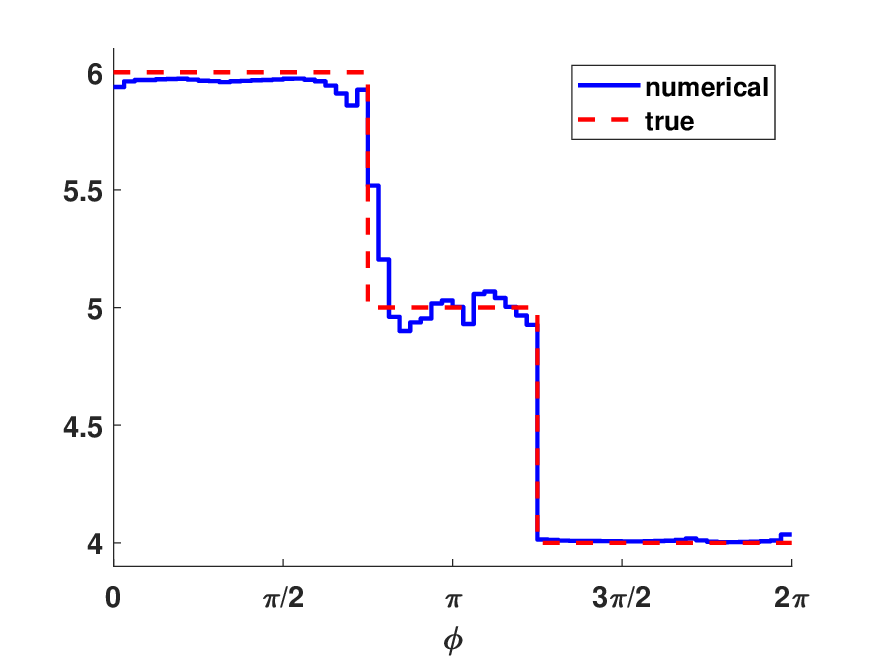}
\caption{Coefficient $a_{11}$ at the $128$-th iteration in the $\phi$-direction of Example \ref{example1}}\label{fig::a11_1}
\endminipage
\end{figure}

The measurement noise model is $z|_{\Xi} = u[a^{\star}]|_{\Xi} + \delta \xi_U || u[a^{\star}] ||_{\infty}$, where $\xi_U \sim \text{Uniform}([-1,1]^{M_1})$ and $\delta=0.001$, preserving $0.1\%$ relative noise level in the infinity norm. For $m=64$,
we initialize all coefficients uniformly $a_1^0=...=a_{64}^0=\left(
\begin{matrix}
5 & 1 & 1\\
1 & 4 & 1\\
1 & 1 & 3
\end{matrix}\right)$. After $\tilde{k}=128$ iterations, Figure \ref{fig::J_1} shows the objective function $J(k)$ evolution, Figure \ref{fig::a11_1} compares reconstructed $a_{11}^{128}$(blue) vs true $a_{11}^{\star}$(red dashed), while other components are omitted. Finally we achieve relative error $\|a_h^{128} - a^{\star}\|_{1}/\|a^{\star}\|_{1} = \text{7.201e-3}$.

\begin{figure}[H]
\centering
\includegraphics[scale=0.4]{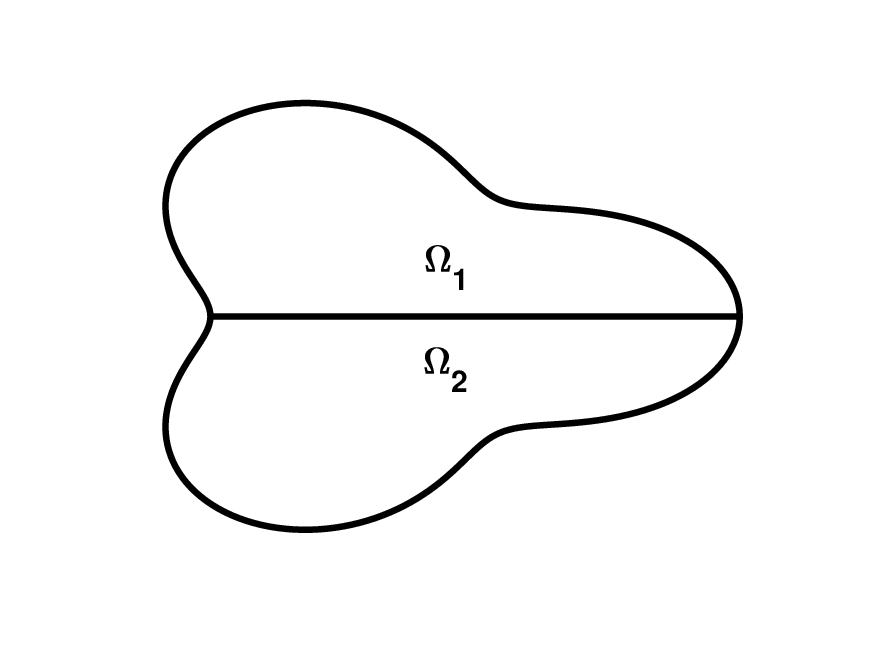}
\caption{Domain $\Omega$ in Example \ref{example4}} \label{fig:domainExample4}
\end{figure}

\begin{example}\label{example4}
The impact of noise in measurement data on our method's performance deserves focused investigation. Assign $\tilde{r}(\phi) = \sqrt{2+ \cos(3\phi)}$ with interfacial angle $\theta_1=0$, $\theta_2=\pi$; see Figure \ref{fig:domainExample4}. Set $f=(1, 1)^T$, $p =(1, 1)^T$ and $a_1^{\star}=\left(
\begin{matrix}
8 & 1 & 1\\
1 & 6 & 1\\
1 & 1 & 4
\end{matrix}\right)$,$a_2^{\star}=\left(
\begin{matrix}
4 & 1 & 0.5\\
1 & 3 & 0.5\\
0.5 & 0.5 & 2
\end{matrix}\right)$.
\end{example}

The measurement noise model is $z|_{\Xi} = u[a^{\star}]|_{\Xi} + \delta \xi_U || u[a^{\star}] ||_{\infty}$, where $\xi_U \sim \text{Uniform}([-1,1]^{M_1})$ and $\delta=0.001,0.003,0.005$. For each noise level $\delta$, we set $m=16$ and initialize the algorithm from $a_1^{0}=...=a_{16}^0=\left(
\begin{matrix}
6 & 0.75 & 0.75\\
0.75 & 4.5 & 0.75\\
0.75 & 0.75 & 3
\end{matrix}\right)$. Table \ref{tab:case1} reports the number of iterations $\tilde{k}$ required before termination and the $L^1$ relative errors in the anisotropic coefficients at $\tilde{k}$-th iteration for each $\delta$. 
\begin{table}
\centering
\caption{Iteration $\tilde{k}$ and $L^1$ relative errors of Example \ref{example4} for different $\delta$}
\label{tab:case1}
\footnotesize
\begin{tabular}{@{}lll}
\br
$\delta$ & $\tilde{k}$ & $\|a_h^{\tilde{k}} - a^{\star}\|_{1}/ \|a^{\star}\|_{1}$ \\
\mr
0.001 &  725 &  3.602e-2 \\
0.003 & 871 & 3.220e-2 \\
0.005 &  755 & 3.792e-2 \\
\br
\end{tabular}
\end{table}

Across Examples \ref{example1}–\ref{example4}, our method consistently recovers both interfaces and  the six anisotropic coefficients, with Example \ref{example4} further illustrating its stability under noisy conditions.

\section{Conclusion} \label{sec:conclusion}
In this paper, we present a method for identifying heterogeneous and linearly anisotropic elastic parameters using only a single full-field measurement. Our approach combines the semi-analytical forward solver with an optimization framework based on a dedicated energy functional incorporating least-squares fitting and total variation regularization, solved efficiently via the Adam algorithm. A key advantage of our method is its computational efficiency: each iteration requires just one forward problem solution and completely avoids adjoint computations.

To address challenges posed by anisotropy, heterogeneity, irregular domains, and solution singularities, we first parameterize the elastic domain boundary as a (piecewise) $C^1$ mapping. Then establish a curvilinear coordinate system in which the irregular stellar domain is diffeomorphically transformed into a standardized semi-infinite rectangular region. And yield a hybrid formulation combining variational principles with differential operators for linear elasticity. Through angular discretization, the resulting semi-discrete formulation is mathematically tractable via our direct computational approach. Numerical experiments demonstrate that this approach achieves rapid convergence of discretized eigenvalues to their continuous counterparts in the elliptic operator's spectrum while naturally resolving solution singularities. The method also provides optimal error estimates and shows excellent accuracy in forward problem solutions.

Our numerical results confirm the method's effectiveness in simultaneously recovering both interface locations and all six anisotropic coefficient values across subdomains. Future work will extend this approach to linear elastic problems featuring multiple singularities in more general domains.

\ack
This work was partially supported by the NSFC Projects No. 12025104 and the Liaoning Provincial Department of Education Project No. LJ212410147018.

\end{document}